\documentclass[12pt,reqno,a4wide]{amsart}

\usepackage{tikz} 
\usepackage{calrsfs}
\usepackage[charter]{mathdesign}

\allowdisplaybreaks

\title{On third-order Pell polynomials}

\author{Helmut Prodinger}
\address{Helmut Prodinger\\
	Department of Mathematical Sciences\\
	Stellenbosch University\\
	7602 Stellenbosch\\
	South Africa}
\email{hproding@sun.ac.za}


\begin{document}
\begin{abstract}
Binet formulae for three versions of third-order Pell polynomials are derived.
	\end{abstract}

\maketitle
	
\section{A Binet formula for third-order Pell polynomials}

Mahon and Horadam \cite{horadam-ternary} study the recursion for $r_n=r_n(x)$
\begin{equation*}
r_{n}=2xr_{n-1}+r_{n-3},\ n\ge3,\	\quad r_0=0,\ r_1=1,\ r_2=2x.
\end{equation*}
The characteristic equation of the recursion is
\begin{equation*}
X^3-2xX^2-1=0;
\end{equation*}
in \cite{horadam-ternary} the authors try to attack it with Cardano's formula (and without using computers).
Based on our recently acquired experience with ternary recursions, see, e. g., \cite{Rocky}, and using \textsf{Maple},
we are able to get stronger and quite satisfying results.

First,	substitute $X=x/Y$
	and the equation is now
\begin{equation*}
Y^3+2x^3Y-x^3=0.
\end{equation*}
Now, setting $x^3=-1/z$ it is
	\begin{equation*}
	zY^3-2Y+1=0
	\end{equation*}
Using the substitution
\begin{equation*}
z=(1-t)^2(1+t),
\end{equation*}
the equation has beautiful roots:
\begin{align*}
v_1&:=\frac1{1-t},\\
v_2&:=\frac{1+t+\sqrt{(1+t)(5-3t)}}{2(t^2-1)},\\
v_3&:=\frac{1+t-\sqrt{(1+t)(5-3t)}}{2(t^2-1)}.
\end{align*}
Note that $v_2+v_3=1/(t-1)$ and $v_2v_3=1/(t^2-1)$.
 Setting
	\begin{equation*}
r_n=x^{n-1}\Big(Av_1^{-n}+Bv_2^{-n}+Cv_3^{-n}\Big).
	\end{equation*}
we can now compute the coefficients $A,B,C$ from the initial values, 
	with the result
	\begin{align*}
A&=-\frac1{1+3t},\\*
B&=\frac1{2(1+3t)}-\frac32{\frac {\sqrt { ( 1+t  )   ( 5-3t) }}{
		\left( 5-3t \right)  \left( 1+3t \right) }},\\*
C&=\frac1{2(1+3t)}+\frac32{\frac {\sqrt { ( 1+t  )   ( 5-3t) }}{
		\left( 5-3t \right)  \left( 1+3t \right) }}.
	\end{align*}
	Even more appealing are the reciprocal roots:
	\begin{align*}
w_1=v_1^{-1}&=1-t\\
w_2=v_2^{-1}&=\frac{1+t-\sqrt{(1+t)(5-3t)}}{2}\\
w_3=v_3^{-1}&=\frac{1+t+\sqrt{(1+t)(5-3t)}}{2}\\
	\end{align*}
	The promised Binet formula is
	\begin{equation*}
r_n=x^{n-1}\Big(Aw_1^{n}+Bw_2^{n}+Cw_3^{n}\Big).
	\end{equation*}
As an example, 
\begin{equation*}
r_{18}=131072{x}^{17}+245760{x}^{14}+159744{x}^{11}+42240{x}^{8}+4032
{x}^{5}+84{x}^{2}
\end{equation*}
It starts with the highest power $x^{17}$, and goes down from there by powers of $x^{-3}$, as predicted.

\section{The coefficients}

First, we want to show, that as in the Binet formula for Fibonacci numbers, the square root is superficial
and cancels out. To see this, we compute
	\begin{align*}
\frac{(1+3t)r_n}{x^{n-1}}&=-(1-t)^n\\
&+\bigg[\frac1{2}-\frac32{\frac {\sqrt { ( 1+t  )   ( 5-3t) }}{
		\left( 5-3t \right)   }}\bigg]\bigg(\frac{1+t-\sqrt{(1+t)(5-3t)}}{2}\bigg)^n\\
	&+\bigg[\frac1{2}+\frac32{\frac {\sqrt { ( 1+t  )   ( 5-3t) }}{
			\left( 5-3t \right)   }}\bigg]\bigg(\frac{1+t+\sqrt{(1+t)(5-3t)}}{2}\bigg)^n.
	\end{align*}
	The non-trivial part (lines 2 and 3), multiplied by $2^{n+1}$ and the abbreviation $W=\sqrt{(1+t)(5-3t)}$ is
	\begin{align*}
	\Xi:=&\bigg[1-{\frac {3W}{			  5-3t    }}\bigg]\big(1+t-W\big)^n+ \bigg[1+{\frac {3W}{			  5-3t    }}\bigg]\big(1+t+W\big)^n\\
	&=2\sum_k \binom{n}{2k}W^{2k}(1+t)^{n-2k}
	+2\sum_k \binom{n}{2k+1}W^{2k+1}(1+t)^{n-1-2k}\frac {3W}{5-3t}\\
	&=2\sum_k \binom{n}{2k}(5-3t)^k(1+t)^{n-k}
+6\sum_k \binom{n}{2k+1}(5-3t)^{k}(1+t)^{n-k},
	\end{align*}
	and no more square roots are present.
	
	Now we want to show how to compute various series expansions. We start with
	\begin{equation*}
z=(1-t)^2(1+t).
	\end{equation*}
So, $z$ is given in terms of $t$, but we want the expansion of $t$ in terms of $z$.
This can be seen in the context of the Lagrange inversion or the Lagrange-B\"urmann formula, \cite{Henrici-Vol1}.
We will use contour integration to get the inverted series. We need another substition:
	\begin{equation*}
t=u-1\quad \Longrightarrow  \quad z=u(u-2)^2.
	\end{equation*}
\begin{figure}[h]
\includegraphics[width=4cm,height=4cm ]{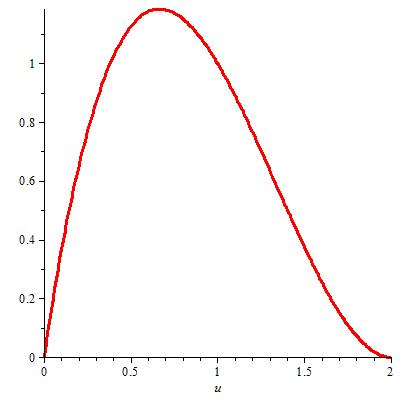}	
\caption{$z$ as a function of $u$. If $|u|<\frac23$ the function can be inverted.}
\end{figure}
The advantage is that $z\to0 \Leftrightarrow u\to0$; we restrict ourselves to $|u|<\frac23$. Furthermore $\dfrac{dz}{du}=(3u-2)(u-2)$. We are going to show that
	\begin{equation*}
u=\sum_{n\ge1}\frac1{n}\binom{3n-2}{n-1}\frac1{2^{3n-1}}z^n.
	\end{equation*}
The usual convergence test shows that this series converges for $|z|<\frac{32}{27}\approx1.185$.

	For that, we compute (the path of integration is in all instances a small circle around the origin)
	\begin{align*}
[z^n]u&=\frac1{2\pi i}\oint \frac{dz}{z^{n+1}}u\\
&=\frac1{2\pi i}\oint \frac{du(3u-2)(u-2)}{u^{n+1}(u-2)^{2n+2}}u\\
&=[u^{n-1}]\frac{3u-2}{(u-2)^{2n+1}}
=\frac1{2^{2n+1}}[u^{n-1}]\frac{2-3u}{(1-\frac u2)^{2n+1}}.
	\end{align*}
	Reading off this coefficient and simplifying leads to the result. Now we will do a similar computation to show that
		\begin{equation*}
-\frac{(1-t)^n}{1+3t}=\sum_{\ell\ge0}2^{n-3\ell-1}\binom{3\ell-n}{\ell}z^\ell,
	\end{equation*}
which is the first term of the Binet formula. First, 
$$
-\frac{(1-t)^n}{1+3t}=-\frac{(2-u)^n}{3u-2},
$$
 and then
\begin{align*}
-[z^\ell ]\frac{(2-u)^n}{3u-2}&=-\frac1{2\pi i}\oint \frac{dz}{z^{\ell+1}}\frac{(2-u)^n}{3u-2}
=-\frac1{2\pi i}\oint \frac{du(3u-2)(u-2)}{u^{\ell+1}(u-2)^{2\ell+2}}\frac{(2-u)^n}{3u-2}\\
&=-(-1)^n\frac1{2\pi i}\oint \frac{du}{u^{\ell+1}(u-2)^{2\ell+1-n}}
=[u^\ell]\frac{1}{ (2-u)^{2\ell+1-n}}\\
&=\frac1{2^{3\ell+1-n }}\binom{3\ell-n}{\ell}.
\end{align*}
If only this first term of the Binet formula would be used, we would have derived that
	\begin{equation*}
	r_n(x)=\sum_{\ell \ge0} \binom{n-1-2\ell }{\ell}(2x)^{n-1-3\ell}.
\end{equation*}
However, $r_n(x)$ isn't an infinite series, it is just a polynomial. This means that the second and third term in the Binet
formula kill off the infinite rest of the series. 
This can be seen for instance, by finding the expansions of $w_2^n+w_3^n$ and $\frac{w_2^n-w_3^n}{w_2-w_3}$ using the Girard-Waring formula, viz.~\cite{Gould}.

Instead of showing this cancellation directly, we prove that
	\begin{equation*}
r_n(x)=\sum_{0\le \ell \le (n-1)/3} \binom{n-1-2\ell }{\ell}(2x)^{n-1-3\ell}
	\end{equation*}
	by simple induction.
	The initial conditions match, and then:
	\begin{align*}
r_{n}&=2xr_{n-1}+r_{n-3}\\
&=2x\sum_{0\le \ell \le (n-2)/3} \binom{n-2-2\ell }{\ell}(2x)^{n-2-3\ell}+\sum_{0\le \ell \le (n-4)/3} \binom{n-4-2\ell }{\ell}(2x)^{n-4-3\ell}\\
&= \sum_{0\le \ell \le (n-2)/3} \binom{n-2-2\ell }{\ell}(2x)^{n-1-3\ell}+\sum_{1\le \ell \le (n-1)/3} \binom{n-2-2\ell }{\ell-1}(2x)^{n-1-3\ell}\\
&= (2x)^{n-1}+\sum_{1\le \ell \le (n-2)/3} \binom{n-1-2\ell }{\ell}(2x)^{n-1-3\ell}+
[\![3\mid(n-1)]\!]\\
&=\sum_{0\le \ell \le (n-1)/3} \binom{n-1-2\ell }{\ell}(2x)^{n-1-3\ell},
	\end{align*}
	as predicted.

	\section{The second version of the third-order Pell polynomials}
	
	We consider the second version proposed by  Mahon and Horadam, but only indicate what changes.
	These polynomials are given by
	\begin{equation*}
		s_{n}=2xs_{n-1}+s_{n-3},\ n\ge3,\	\quad r_0=0,\ r_1=2,\ r_2=2x.
	\end{equation*}
	Since the recursion is the same, only the initial values change, and they are
	\begin{align*}
		A&=\frac{2t}{(1+3t)(t-1)}\\*
		B&=\frac t{(1+3t)(t-1)}+ \frac{3t^2-3t-2}{(t^2-1)(3t+1)(3t-5)}\\*
		C&=\frac t{(1+3t)(t-1)}- \frac{3t^2-3t-2}{(t^2-1)(3t+1)(3t-5)}
	\end{align*}
	The Binet formula is then
		\begin{equation*}
		s_n=x^{n-1}\big(Aw_1^{n}+Bw_2^{n}+Cw_3^{n}\big).
	\end{equation*}
	There is also an explicit formula, which works for $n\ge2$:
	\begin{equation*}
s_n(x)=\sum_{0\le \ell \le (n-1)/3} (2x)^{n-1-3\ell}\frac{n-\ell-1}{n-2\ell-1}\binom{n-2\ell-1}{\ell}
	\end{equation*}
	As before, this could be guessed from the first term in the Binet formula and proved by induction.
	
	\section{The third version of the third-order Pell polynomials}
	
These polynomials  are given by
\begin{equation*}
	\sigma_{n}=2x\sigma_{n-1}+\sigma_{n-3},\ n\ge3,\	\quad \sigma_0=3,\ \sigma_1=2x,\ \sigma_2=4x^2.
\end{equation*}
This time the Binet formula is very simple:
\begin{equation*}
\sigma_n=x^n(w_1^n+w_2^n+w_3^n).
\end{equation*}
	There is an explicit formula, valid for $n\ge1$:
	\begin{equation*}
\sigma_n=
\sum_{0\le \ell \le n/3} \frac{n}{n-2\ell}\binom{n-2\ell}{\ell}(2x)^{n-3\ell}
	\end{equation*}
	\bibliographystyle{plain}

\end{document}